\documentclass[11pt]{amsart}

\usepackage{amssymb, amsmath, amscd, amsthm, graphicx, psfrag}
\usepackage[matrix,arrow,curve]{xy}
\xyoption{dvips}              

\newcommand\ie{{\it i.e. }}

\newcommand{\R}{{\mathbb R}}

\newcommand{\Q}{{\mathbb Q}}
\newcommand{\calL}{{\mathcal{L}}}

\newcommand{\calP}{{\mathcal{P}}}
\newcommand{\calE}{{\mathcal{E}}}
\def\nbox{\hfill $\Box$}

\newcommand\rth{\refstepcounter{equation}}
\newcommand\numb{\rth{\rm \theequation}}
\numberwithin{equation}{section}

\DeclareMathOperator{\Tot}{Tot}

\renewcommand{\Im}{\mathop{\mathrm{Im}}\nolimits}

\newcommand{\ASS}{{\mathcal{ASS}}}
\newcommand{\COASS}{{\mathcal{COASS}}}
\newcommand{\unit}{{\mathbf{1}}}

\newcommand{\Vect}{\mathrm{Vect}}

\newcommand{\CDGA}{\mathrm{CDGA}}

\newcommand{\calC}{{\mathcal{C}}}

\newcommand{\calA}{{\mathcal{A}}}
\newcommand{\calD}{{\mathcal{D}}}
\newcommand{\calO}{{\mathcal{O}}}

\theoremstyle{plain}
\newtheorem{thm}{Theorem}[section]
\newtheorem{theorem}[thm]{Theorem}

\newtheorem{proposition}[thm]{Proposition}
\newtheorem{lemma}[thm]{Lemma}

\newtheorem{lemma-def}[thm]{Lemma-definition}

\theoremstyle{definition}
\newtheorem{definition}[thm]{Definition}
\newtheorem{remark}[thm]{Remark}

\begin{document}


\title{Homotopy graph-complex for configuration and knot spaces}


\author{Pascal Lambrechts}
\address{Institut Math\'{e}matique, 2 Chemin du Cyclotron, B-1348 Louvain-la-Neuve, Belgium}
\email{lambrechts@math.ucl.ac.be}
\urladdr{http://milnor.math.ucl.ac.be/plwiki}
\author{Victor Turchin}
\thanks{The second author was supported in part by the grants NSH-1972.2003.01, RFBR 05-01-01012a.}
\address{University of Oregon, USA. Institut des Hautes Etudes Scientifiques, France.}
\email{vitia-t@yandex.ru} \urladdr{http://www.math.ucl.ac.be/membres/turchin/}
\subjclass{Primary: 57Q45; Secondary: 55P62, 57R40}
\keywords{knot spaces, embedding calculus, Bousfield-Kan spectral sequence,
graph-complexes}


\begin{abstract}
 In the paper we prove that the primitive part of the Sinha homology spectral
 sequence $E^2$-term for the space of long knots is rationally isomorphic to  
 the homotopy $\calE^2$-term.
 We also define natural graph-complexes computing the
 rational homotopy of configuration and of knot spaces.
\end{abstract}

\maketitle



\sloppy

\section{Introduction}\label{s1}

In~\cite{Sinha-TKS,Sinha-OKS} D.~Sinha defined a cosimplicial model for the space $Emb$ of long knots $\R\hookrightarrow\R^d$,
$d\geq 4$. It was proven in~\cite{ALV,LTV} that the associated homology Bousfield-Kan spectral sequence collapses rationally
at the second term. The same result was established for the associated homotopy spectral sequence (over $\Q$). The proof
will appear in~\cite{ALTV-coformal}. But $Emb$ is an $H$-space with a homotopy commutative product\footnote{It was shown recently by
P.~Salvatore that $Emb$ is a double loop space~\cite{Salvatore}.}. It implies in particular that 
the $\calE^2$ term of the
(co)homotopy spectral sequence must be rationally isomorphic to the primitive part of the $E^2$ (co)homology term.

\begin{figure}[h!]
\psfrag{p}[0][0][1][0]{$p$}
\psfrag{q}[0][0][1][0]{$q$}
\psfrag{p'}[0][0][1][0]{$p'$}
\psfrag{q'}[0][0][1][0]{$q'$}
\includegraphics[width=8cm]{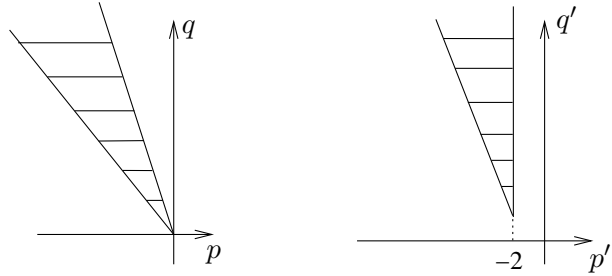}
\caption{Homology $E^2_{p,q}$ and homotopy $\calE^2_{p',q'}$ terms.}\label{fig1}
\end{figure}
 
The (co)homology $E^2_{p,q}$-term is concentrated in the second quadrant between two 
lines~\cite{T-OS}:
$$
\begin{array}{rll}
q=&-\frac{d-1}2 p& \qquad \text{(lower line)}\\
q=&-(d-1)(p+1)& \qquad \text{(upper line)}
\end{array}
$$
The (co)homotopy $\calE^2_{p',q'}$-term is also concentrated in the second quadrant and bounded by the lines~\cite{SS}:
$$
\begin{array}{rll}
q'=&-(d-2)p'-d+3& \quad \text{(lower line)}\\
p'=&-2& \quad \text{(right line)}
\end{array}
$$


In Part~I of the paper we will give a simple and purely algebraic proof of this isomorphism:
$$
\calE^2_{*,*}=Prim(E^2_{*,*}).
$$
In particular we will see in Section~\ref{s4} how via this 
isomorphism the bigradings of both spectral sequences are related to each other. We will see that the lower line of the homotopy spectral sequence corresponds to the lower line of the homology one.
The right line $p'=-2$ of $\calE^2$ corresponds to the upper line $q=-(d-1)(p+1)$ of $E^2$. In general any vertical line $p'=-n$
corresponds to $q=-(d-1)(p+n-1)$.

This isomorphism in the case of the lower lines (on the level of the bialgebra of chord diagrams) was proved by
J.~Conant~\cite{Conant}. He gives an elegant reformulation of his result using 3-valent graphs.

Part~II, which actually gave the name to the paper, is devoted to graph-complexes. Our motivation was to produce new
graph-complexes whose homology has a nice geometrical interpretation. We define a series of graph-complexes that compute the
rational homotopy of configuration spaces. Building up on this series of complexes we define a bigger complex whose homology is
the rational homotopy of the space of long knots. A more thorough introduction for Part~II is Section~\ref{s7}.

\part{Isomorphism $\calE^2=Prim(E^2)$}

\section{Cosimplicial model for the space of long knots modulo immersions}\label{s2}

The space $\overline{Emb}$ of long knots modulo immersions is the homotopy fiber of the inclusion
$$
Emb\hookrightarrow Imm
$$
of the space of long knots $Emb$ in the space of long immersions $Imm$. By the word \lq\lq long" we understand smooth map
$\R^1\to\R^d$ that coincide with a fixed linear map outside a compact subset of $\R^1$. (We will deliberately omit $d$ 
to simplify notation, assuming that the dimension $d\geq 4$ is fixed once and forever.\footnote{When we consider 
configuration spaces we assume $d\geq 3$.})

D.~Sinha showed in~\cite{Sinha-OKS} that $\overline{Emb}$ is homotopy equivalent to $Emb\times \Omega Imm\simeq Emb\times 
\Omega^2 S^{d-1}$. So, the homology and homotopy of $Emb$ are easily related to those of $\overline{Emb}$ and the results for $\overline{Emb}$ that we obtain in the first part of the paper can be obviously reestablished for $Emb$. 

In~\cite{Sinha-TKS} D.~Sinha defined a cosimplicial space whose homotopy totalization is $\overline{Emb}$. The $n$-th component
$C^n$ of the cosimplicial space is some compactification of the configuration space of points in
$I\times\R^{d-1}=[0,1]\times\R^{d-1}$:
$$
\left\{ (x_0,x_1,\ldots,x_{n+1})\left| 
\begin{array}{ll}
x_i\in I\times\R^{d-1}; & x_i\neq x_j\\
x_0=(0,\bar 0);& x_{n+1}=(1,\bar 0)
\end{array}
\right.
\right\}
$$

Given $0\leq i\leq n+1$, the coface map $d^i$ is doubling of the $i$-th point $x_i$ of the configuration in the 
direction $(1,\bar0)$. The codegeneracy map $s^i$, for $i=1\ldots n$, is given by forgetting of $x_i$.

In~\cite{Sinha-OKS} D.~Sinha provides another \lq\lq operadic" construction for the cosimplicial replacement of $\overline{Emb}$
(this cosimplicial space is homotopy equivalent to the previous one).

For any symmetric monoidal category $(\calC,\otimes,\unit)$ denote by $\ASS=\{\unit\}_{n\geq 0}$ the 
associative non-$\Sigma$ operad.\footnote{Dually we will denote by $\COASS=\{\unit\}_{n\geq 0}$ the 
associative non-$\Sigma$ cooperad in $(\calC,\otimes,\unit)$.}

Provided an operad $\calO$ in $(\calC,\otimes,\unit)$ is endowed with a morphism $\ASS\to\calO$,
the collection $\{\calO^n\}_{n\geq 0}=\{\calO(n)\}_{n\geq 0}$ becomes a cosimplicial object in this category.
Cofaces $d^i\colon\calO^n\to\calO^{n+1}$, $i=0\ldots n+1$, are compositions with $m=\unit=\ASS(2)$:
$$
d^0({-})=m\circ_2{-}; \quad d^i({-})={-}\circ_i m,\,\, i=1\ldots n; \quad d^{n+1}({-})=m\circ_1{-}.
$$
Codegeneracies are compositions with $e=\unit=\ASS(0)$
$$
s^i({-})={-}\circ_i e,\,\, i=1\ldots n.
$$

Sinha applies this standard construction to an operad $\{C\langle n\rangle\}_{n\geq 0}$ homotopy equivalent to the operad of
little $d$-cubes. Each space $C^n=C\langle n\rangle$ of this operad is the compactification in
 $(S^{d-1})^{n\choose 2}$
of the space of reciprocal directions of $n$ distinct points in~$\R^d$:
$$
\left\{ \left(\frac{x_j-x_i}{|x_j-x_i|}\right)_{1\leq i<j\leq n} \left| 
\begin{array}{c}
x_i\in\R^d;\\
x_i\neq x_j
\end{array}
\right.
\right\}
\subset (S^{d-1})^{n\choose 2}.
$$

We will assume that we work with one of these cosimplicial models. The operadic and cosimplicial structures of $C^\bullet$
induce similar structures on the (co)homology and (co)homotopy of $C^\bullet$. The cohomology simplicial algebra
will be denoted by
$$
A_\bullet=\{A_n\}_{n\geq 0}=\{H^*(C^n)\}_{n\geq 0}.
$$
We will consider the cohomology only with rational coefficients.

The rational cohomotopy simplicial Lie coalgebra will be denoted by
$$
L_\bullet=\{L_n\}_{n\geq 0}=\{Mor(\pi_*(C^n),\Q)\}_{n\geq 0}.
$$

We will also work with the rational homotopy cosimplicial Lie algebra
$$
L^\bullet=\{L^n\}_{n\geq 0}=\{\pi_*(C^n)\otimes\Q\}_{n\geq 0}.
$$

\section{Explicit description of $A_n=H^*(C^n)$ and of $L^n=\pi_*(C^n)\otimes\Q$}\label{s3}

The algebras $A_n$, $n\geq 0$, are well known~\cite{Arn,Coh}. Being graded commutative they are generated by
$a_{ij}$, $1\leq i\neq j\leq n$, of degree $d-1$, that satisfy the relations:
$$
\begin{array}{cl}
& \,\quad a_{ji}=(-1)^d a_{ij}\\
\begin{array}{c}
\text{quadratic}\\
\text{relations}
\end{array}
&
\left\{
\begin{array}{l}
a_{ij}^2=0\\
a_{ij}a_{jk}+a_{jk}a_{ki}+a_{ki}a_{ij}=0
\end{array}
\right.
\end{array}
\eqno(\numb)\label{eq31}
$$

We assume that the component $C\langle 0\rangle$ is a point, so $A_0=\Q$.

\begin{figure}[h!]
\includegraphics[width=5cm]{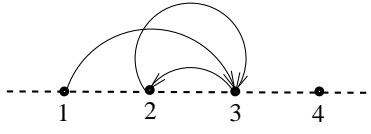}
\caption{Graph corresponding to $a_{13}a_{32}a_{23}\in A_4$.}\label{fig2}
\end{figure}

Any monomial can be viewed as a directed graph on the set $\{1,2,\ldots,n\}$: the directed edge $(i,j)$ is put exactly the number
of times the generator $a_{ij}$ is represented in the monomial, see Figure~\ref{fig2}.

It can be easily seen that a monomial in $A_n$ is non-zero if and only if the corresponding graph is a forest. On such a 
graph $\Gamma$ the face map $d_0(\Gamma)$ is non-zero if and only if the valence of vertex~1 is zero. In this case $d_0$ simply
removes the vertex~1, all other vertices are shifted by~1. Face $d_i(\Gamma)$, $i=1\ldots n-1$, is obtained by collapsing the
segment $[i,i+1]$. If $\Gamma$ contains the edge $(i,i+1)$, then $d_i(\Gamma)=0$. Face $d_n$ acts similarly to $d_0$ removing the
last $n$-th vertex.

\vspace{0.3cm}

\begin{figure}[h]
\psfrag{d0}[0][0][1][0]{$d_0$}
\psfrag{d1}[0][0][1][0]{$d_1$}
\psfrag{d2}[0][0][1][0]{$d_2$}
\psfrag{d3}[0][0][1][0]{$d_3$}
\includegraphics[width=13cm]{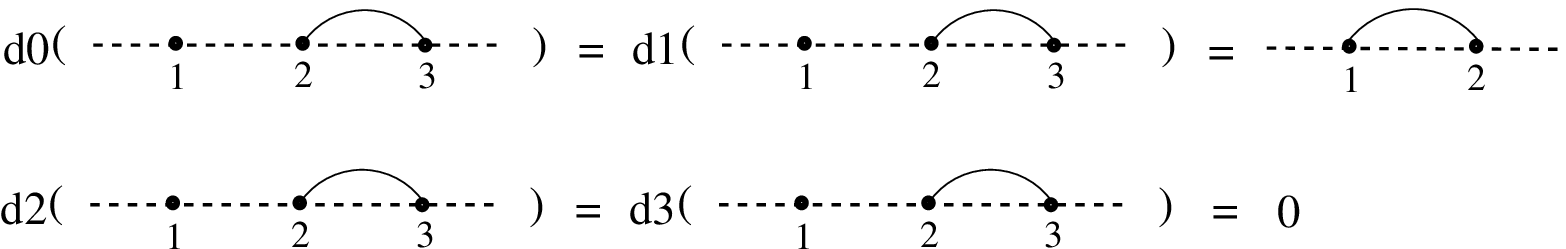}
\caption{}\label{fig3}
\end{figure}


The degeneracy map $s_i$, for $i=1\ldots n+1$, inserts a new vertex between $i-1$ and $i$.

\vspace{0.3cm}

\begin{figure}[h]
\psfrag{s1}[0][0][1][0]{$s_1$}
\psfrag{s2}[0][0][1][0]{$s_2$}
\psfrag{s3}[0][0][1][0]{$s_3$}
\includegraphics[width=15cm]{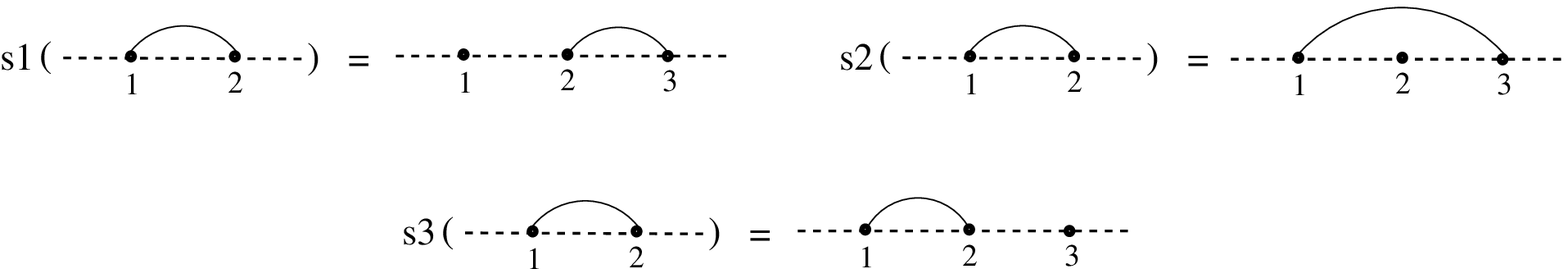}
\caption{}\label{fig4}
\end{figure}


The normalized part $NA_n$ of $A_n$
$$
NA_n=A_n/+_{i=1}^n\Im s^i
$$
is spanned by the forests with each vertex $1,\ldots,n$ of positive valence.

The first term $E_1=(\oplus_{p=0}^\infty E_1^{-p,*},d)$ of the cohomology Sinha spectral sequence is the normalized
complex $\Tot A_\bullet = (\oplus_{p=0}^\infty s^{-p}NA_p,d)$, where $s^{-p}$ denotes $p$-fold desuspension. The differential is as 
usual the alternated sum of faces $d_i$.

\vspace{0.5cm}

The Lie algebra $L^n$, $n\geq 0$, is generated by $\alpha_{ij}$, 
$1\leq i\neq j\leq n$, of degree $d-1$~\cite{Kohno}. The relations are 
$$
\begin{array}{cl}
& \,\,\quad\alpha_{ji}=(-1)^d \alpha_{ij}\\
\begin{array}{c}
\text{quadratic}\\
\text{relations}
\end{array}
&
\left\{
 \begin{array}{l}
 \left[ \alpha_{ij},\alpha_{kl}\right]=0, \quad\text{ if  $\#\{i,j,k,l\}=4$}\\
 \left[ \alpha_{ij},\alpha_{jk}+\alpha_{ki}\right]=0
 \end{array}
\right.
\end{array}
\eqno(\numb)\label{eq32}
$$

The bracket in $L^n$ is the Whitehead bracket which is of degree $-1$.

\vspace{0.5cm}

It is well known that $A_n$ and $L^n$ are Koszul dual~\cite{CohG}. This means that the ${n\choose 2}$-dimensional space $V_n$ of generators of $A_n$ 
 is dual to the space $V^n$  of generators of $L^n$. And the space $R_n\subset S^2V_n$ 
 spanned 
by the quadratic relations~\eqref{eq31} is orthogonal to the space $R^n\subset S^2V^n$ of quadratic relations of $L^n$.
Moreover $A_n$ and $L^n$ are Koszul which means some nice homological property of their bar-constructions.
This property will be used in the proof of our main result Theorem~\ref{t31}.

\vspace{0.5cm}

The cofaces $d^k\colon L^n\to L^{n+1}$, for $k=0\ldots n+1$; and the codegeneracies $s^k\colon L^n\to L^{n-1}$, for $k=1\ldots n$,
are defined on generators as follows:
$$
\begin{array}{ll}
d^k(\alpha_{ij})=
\begin{cases}
\alpha_{ij},&\text{if $i<j<k$;}\\
\alpha_{ij}+\alpha_{i,j+1},&\text{if $i<j=k$;}\\
\alpha_{i,j+1},&\text{if $i<k<j$;}\\
\alpha_{i,j+1}+\alpha_{i+1,j+1},&\text{if $i=k<j$;}\\
\alpha_{i+1,j+1},&\text{if $k<i<j$.}
\end{cases}
&
\quad s^k(\alpha_{ij})=
\begin{cases}
\alpha_{ij},&\text{if $i<j<k$;}\\
0,&\text{if $i<j=k$;}\\
\alpha_{i,j-1},&\text{if $i<k<j$;}\\
0,&\text{if $i=k<j$;}\\
\alpha_{i-1,j-1},&\text{if $k<i<j$.}
\end{cases}
\end{array}
$$
In particular $d^0(\alpha_{ij})=\alpha_{i+1,j+1}$, 
$d^{n+1}(\alpha_{ij})=\alpha_{ij}$.

The normalized part
$$
NL^n=\bigcap_{i=1}^n\ker s^i \subset L^n
$$
is spanned by the monomials that use each index $i=1\ldots n$. The space $NL^n$ is isomorphic to a subspace of a graded free Lie algebra 
generated by $x_1=\alpha_{12}$, $x_2=\alpha_{13}$, $\ldots$, $x_{n-1}=\alpha_{1,n-1}$ spanned by the monomials using each $x_i$, $1\leq i\leq n-1$.

The first term $\calE^1=(\oplus_{p=0}^\infty\calE^1_{-p,*},d)$ of the homotopy Sinha spectral sequence is the normalized complex
$
\Tot L^\bullet=(\oplus_{p=0}^\infty s^{-p}NL^p,d).
$

\vspace{0.3cm}

Here is our main result.

\begin{theorem}\label{t31}
(i) The $\calE^2$ term of the homotopy Sinha spectral sequence (for $\overline{Emb}$) is rationally isomorphic to the primitive part
of the homology $E^2$ term.

(ii) The $\calE_2$ term of the cohomotopy Sinha spectral seqence (for $\overline{Emb}$) is rationally isomorphic to the primitive part of the
cohomology $E_2$-term.
\end{theorem}

Since $E^2$-homology term is a polynomial bialgebra\footnote{This is true for any field of coefficients~\cite[Corollary~13.4]{T-one-term}.}, assertions (i) and (ii) 
are equivalent. So, we will prove only (ii). The proof will be given in Section~\ref{s6}.

\section{Correspondence of bigradings}\label{s4}

In this section we describe how the homotopy spectral sequence bigradings $(p',q')$ are related to the homology spectral sequence bigradings $(p,q)$
via the isomorphism of Theorem~\ref{t31}.

We will give an heuristic explanation of this correspondence. But one can easily establish it by a simple analysis of the proof given in
Section~\ref{s6}.

A monomial of degree $i$ in $NA_j$ is an element of $E_1^{p,q}$ with $p=-j$, $q=(d-1)i$. A monomial of degree $i'$ in $NL^{j'}$ is an element
of $\calE^1_{p',q'}$ with $p'=-j'$, $q'=(d-1)i'-(i'-1)=(d-2)i'+1$ (recall that the bracket in $L^\bullet$ is of degree $-1$).

The degree $i$, $i'$ in both cases will be called {\it complexity}. The complexity is preserved by the differential.

Up to a shift of grading the complexes $\Tot A_\bullet$, $\Tot L^\bullet$ depend only on the parity of $d$. So it is natural to expect that
the isomorphism of Theorem~\ref{t31} respects this periodicity and therefore preserves the complexity. The total grading $p+q$, $p'+q'$ must be also unchanged.

Let us find the bigrading $(i',j')=(i,j')$ of $\Tot L^\bullet$ that should correspond to the bigrading $(i,j)$ of $\Tot A_\bullet$. We have
$$
\begin{array}{ccl}
p+q&=&(d-1)i-j\\
p'+q'&=&(d-2)i+1-j'
\end{array}
$$
So $p+q=p'+q'$ implies
$$
j'=j-i+1.
\eqno(\numb)\label{corresp}
$$

For example if $j=2i$ (the case of the lower line in $E_2$, which corresponds to the bialgebra of chord diagrams), one has $j'=i+1$. This corresponds 
to the lower line in $\calE_2$.

The case $j=i+1$ (upper line in $E_2$) produces $j'=2$ (right line in $\calE_2$). This situation produces exactly the homotopy of the factor $\Omega^2S^{d-1}$ of $\overline{Emb}$~\cite{T-OS}. 

In general for a non-trivial monomial of degree $i$ in $NA_j$ the number $j-i$ is the number of connected components in the corresponding forest, see Section~\ref{s3}. So, the correspondence~\ref{corresp} can be resumed as follows: the number of connected components of forests in the cohomological case corresponds to the number of points (of configuration spaces) minus 1 in the homotopy case.

\section{Fixing notations}\label{s5}

In this section we review some of necessary background and fix some notation.

\subsection{$B/B^2$}\label{ss51}

By $\CDGA$ we understand the category of graded connected differential graded algebras with differential raising the degree by~1. Almost all algebras
we deal with are 1-connected, \ie their 1-degree component is trivial. 

Consider a functor from $\CDGA$ to the category of differential graded vector spaces (complexes):
$$
\begin{array}{rccl}
P\colon&\CDGA&\longrightarrow&dg{-}\Vect\\
& B&\longmapsto& B_{>0}/(B_{>0})^2.
\end{array}
$$
For simplicity of notation $P(B)=B_{>0}/(B_{>0})^2$ will be denoted by $B/B^2$.

\subsection{$\calL(B)$}\label{ss52}
By $\calL$ we denote the cobar construction
$$
\calL\colon \CDGA\to dg{-}\mathrm{coLie},
$$
which assigns to any commutative $dg$-algebra $B$ a free $dg$-Lie coalgebra with cobracket of degree~1:
$$
\calL(B)=\bigoplus_{n\geq 1}(\mathrm{coLie}(n)\otimes(B_{>0})^{\otimes n})_{S_n},
$$
whose differential is a sum of two things --- one arising from the initial differential of $B$, the other --- from multiplication in $B$. A nice explicit description of this construction is given in~\cite{Sinha-graphs}.
Notice that in our construction the degree of each space $\mathrm{coLie}(n)$ is $1-n$.

One has a natural transformation
$$
\begin{array}{rccc}
\alpha\colon&\calL&\longrightarrow&P\\
&\calL(B)&\stackrel{\alpha_B}{\longrightarrow}&B/B^2
\end{array}
$$
which is a morphism of complexes sending $\bigoplus_{n\geq 2}(\mathrm{coLie}(n)\times (B_{>0})^{\otimes n})_{S_n}$ to zero and $\mathrm{coLie}(1)\otimes B_{>0}=
B_{>0}$ to the quotient $B_{>0}/(B_{>0})^2$.

The following is a standard result in the rational homotopy theory~\cite{FHT-RHT}.

\begin{proposition}\label{p52}
If $B$ is a polynomial algebra then the map
$\calL(B)\stackrel{\alpha}{\longrightarrow}B/B^2$
is a quasi-isomorphism.
\end{proposition}

\subsection{Totalization}\label{ss53}
Let $V_\bullet$  be a simplicial $dg$-vector space. If not stated otherwise we always assume that the differential raises the degree by~1.
We define $\Tot V_\bullet$ as a complex whose space is $\oplus_{n\geq 0} s^{-n}NV_n$ and the differential is the sum of inner differential of each
$V_n$ plus the alternated sum of faces. Notice that $\Tot V_\bullet$ might be negatively graded, however in all the considered cases the totalization 
always produces positively graded complexes.

Let $V_\bullet$ and $W_\bullet$  be two simplicial $dg$-spaces. Assume that $\Tot V_\bullet$ and  $\Tot W_\bullet$ are left-bounded. One has the 
Eilenberg-MacLane quasi-isomorphism~\cite[\S~29]{May-SO}:
$$
\Tot V_\bullet\otimes \Tot W_\bullet \stackrel{EM}{\longrightarrow} \Tot (V_\bullet\otimes W_\bullet).
$$

This map permits to define a product on the totalization of any simplicial commutative $dg$-algebra $B_\bullet$:
$$
\Tot B_\bullet\otimes \Tot B_\bullet \stackrel{EM}{\longrightarrow} \Tot (B_\bullet\otimes B_\bullet)\stackrel{\mu_\bullet}{\longrightarrow}
\Tot B_\bullet.
$$

The Eilenberg-MacLane map has nice properties. It is associative and $S_n$-equivariant. This implies that $\Tot B_\bullet$ is a commutative
$dg$-algebra. For example, for the simplicial algebra 
$A_\bullet=H^*(C^\bullet)$ the product on $\calA=\Tot A_\bullet$ is the shuffle of diagrams:

\vspace{0.3cm}

\begin{figure}[h]
\psfrag{+}[0][0][1][0]{$\pm$}
\includegraphics[width=15cm]{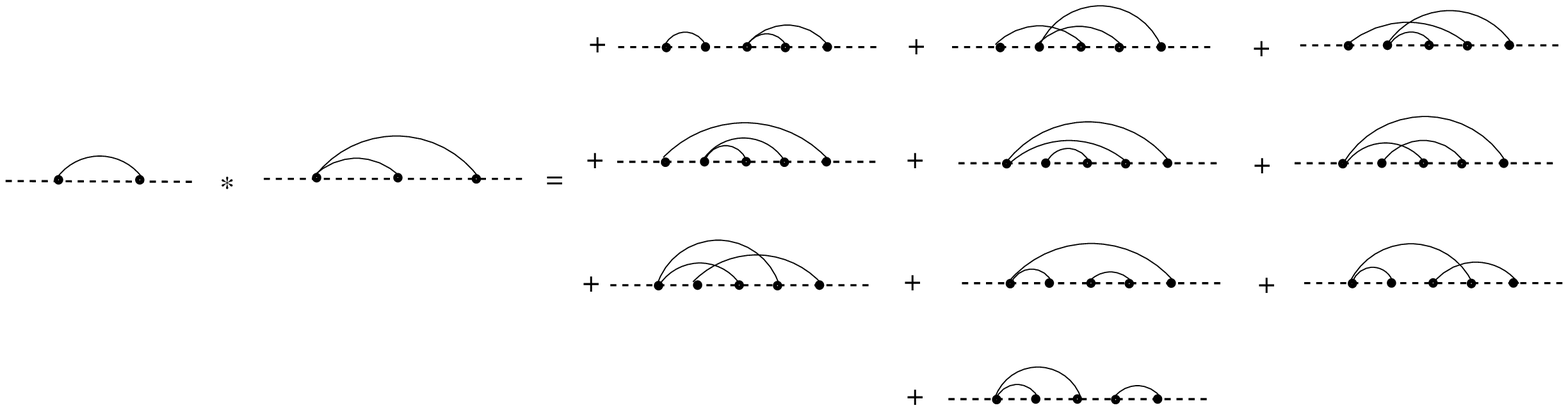}
\caption{The product in $\calA$.}\label{fig5}
\end{figure}


 Consider any polynomial functor
 $$
 \begin{array}{rccl}
 \calP\colon& dg{-}\Vect& \longrightarrow& dg{-}\Vect\\
 & V & \longmapsto& \oplus_{n\geq 0}(\calP(n)\otimes V^{\otimes n})_{S_n}.
 \end{array}
 $$
 
 The property that Eilenberg-MacLane map is associative and $S_n$-equivariant permits to define 
 a morphism:
 $$
 \calP(\Tot V_\bullet)\stackrel{EM_\calP}{\longrightarrow} \Tot \calP(V_\bullet)
 \eqno(\numb)\label{eq51}
 $$
 
 \begin{lemma}\label{l53}
 For a field of characteristic zero, the morphism~\eqref{eq51} is always a quasi-isomorphism 
 for any simplicial $dg$-vector space $V_\bullet$
 (provided $\Tot V_\bullet$ is left-bounded).
 \end{lemma}
 
 \begin{proof}
 One has
 $$
 (\Tot V_\bullet)^{\otimes n}\stackrel{EM_n}{\longrightarrow} \Tot(V_\bullet^{\otimes n})
 $$
 is an $S_n$-equivariant quasi-isomorphism. In characteristic zero it implies that
 $$
 (\calP(n)\otimes(\Tot V_\bullet)^{\otimes n})_{S_n}\stackrel{EM_{\calP(n)}}{\longrightarrow} 
 \Tot((\calP(n)\otimes V_\bullet^{\otimes n})_{S_n})
 $$
 is also a quasi-isomorphism.
 \end{proof}
 
 \begin{lemma}\label{l54}
 For any simplicial commutative $dg$-algebra $B_\bullet$ (provided $\Tot B_\bullet$ is
  positively graded) the map $EM_\calL$ is a quasi-isomorphism:
 $$
 \calL(\Tot B_\bullet)\xrightarrow[EM_\calL]{\simeq} \Tot\calL (B_\bullet).
 $$
 \end{lemma}
 
 \begin{proof}
 First one has to check that $EM_\calL$ is a morphism of complexes. This is so because the 
 product part of the differential in $\calL(\Tot B_\bullet)$
 goes exactly to the product part of the differential in $\Tot\calL(B_\bullet)$ (here one uses 
 the fact that the product in $\Tot B_\bullet$
 was defined through the Eilenberg-MacLane map).
 To see that $EM_\calL$ is isomorphism in homology one can consider  the spectral sequences for 
 both complexes assigned to the filtration by the
 degree $n$ of the polynomial functor $\calL=\oplus_{n\geq 1}\calL_n$. It follows from 
 Lemma~\ref{l53} that the induced map of spectral sequences is
 an isomorphism starting from the first page.
 \end{proof}

 \section{Proof of Theorem~\ref{t31}}\label{s6}
 
 The proof is the following sequence of quasi-isomorphisms.
 $$
 \xymatrix{
 \calA/\calA^2&\calL(\calA)=\calL(\Tot 
 A_\bullet)\ar@{->>}[l]_-\simeq^-\alpha\ar[r]^-\simeq_-{EM_\calL}&
 \Tot\calL(A_\bullet)&\,\Tot(L_\bullet)\ar@{_{(}->}[l]_-\simeq
 }\eqno(\numb)\label{zig-zag}
 $$
 It is well known that $\calA=\Tot A_\bullet$ is a commutative non-cocommutative 
 $dg$-bialgebra, moreover its homology 
 bialgebra 
 $H^*(\calA)$ is polynomial~\cite{T-HLN,T-OS}. 
 We have by Proposition~\ref{p52} that the first arrow $\alpha$ is a quasi-isomorphism.  So the 
 homology of both complexes $\calA$, $\calL(\calA)$ 
 is the space of generators of $H^*(\calA)$ which is exactly the space of primitives.
 
 The second arrow is a quasi-isomorphism by Lemma~\ref{l54}.
 
 Let us explain the last quasi-isomorphism
 $$
 \xymatrix{
 \Tot(L_\bullet)\ar@{^{(}->}[r]^-\simeq& \Tot\calL(A_\bullet).
 }
 \eqno(\numb)\label{eq60}
 $$
The algebra $L^n$  is the Lie Koszul dual of $A_n$~\cite{CohG}. One has the natural inclusion (of Lie coalgebras):
 $$
 L_n\hookrightarrow \calL(A_n).
 \eqno(\numb)\label{eq61}
 $$
 The map~\eqref{eq61} describes the so called \lq\lq diagonal" homology of $\calL(A_n)$. The 
 property $A_n$ is Koszul means that $\calL(A_n)$
 has only diagonal (non-trivial) homology. In other words~\eqref{eq61} is a quasi-isomorphism. 
 But~\eqref{eq61} is a simplicial morphism. As a
 consequence~\eqref{eq60} is also a quasi-isomorphism.
 
 To see that the bigradings correspond by the way described in Section~\ref{s4}, one 
 should generalize the grading {\it complexity} on all the intermediate 
 complexes of the zig-zag~\eqref{zig-zag}, and to show that all the morphismes preserve 
 it. We leave it as an exercise to the reader.
 \nbox

 \part{Graph-complexes}
 
 \section{Introduction}\label{s7}
 Graph-complexes are widely used to study the homology of interesting spaces and to prove 
 interesting theorems~\cite{ConV1,ConV2,CulV,FI,GK,Kont1,Kont2,Kont3,Penner}.
 One series of such graph-comples (its slight modification will be denoted 
 by $\{D_n\}_{n\geq 0}$ 
 throughout the paper) was used by M.~Kontsevich 
  to prove the formality of the operad of little $d$-cubes~\cite{Kon-OMDQ}. The 
 idea of the proof is that $\{D_n\}_{n\geq 0}$ are 
 quasi-isomorphic to the cochains of the operad, and also one has projections inducing homology 
 isomorphism
 $$
 \xymatrix{
 D_n\ar@{->>}[r]^-\simeq_-{\bar I_n}&A_n=H^*(C\langle n\rangle).
 }
 \eqno(\numb)\label{eq71}
 $$
 $A_n$ is considered as a commutative $dg$-algebra with zero differential.
 
 A more thorough account on this result was given by I.~Volic and the first author 
 in~\cite{LambVol-FLBO}.
 
 In~\cite{Catt,Catt2} another graph-complex was defined (its slight modification will be denoted 
 by $\calD$ in the paper). By means of integration
 over configuration spaces this complex was naturally mapped to the De Rahm complex of the space 
 $Emb$ of long knots. One conjectures that this map is a 
 quasi-isomorphism. The reason why it might be so is that $\calD$ is quasi-isomorphic to 
 $\calA=\Tot A_\bullet$, see
 Theorem~\ref{t88}, and therefore the homology of $\calD$ is 
 exactly the rational homology of $Emb$.
 
 The complexes $\{D_n\}_{n\geq 0}$ form a simplicial commutative $dg$-algebra. Its totalization 
 $\Tot D_\bullet$ is exactly the complex $\calD$.
 
 Another motivation for us to study graph-complexes  is that they generalize on  higher 
 homology of knot spaces the 3-valent diagrams calculus developed by Dr.~Bar 
 Natan~\cite{BarNatan} (in the relation with the finite type knot invariants). For 
 example, Theorem~6 of~\cite{BarNatan}, which says that the bialgebra of chord 
 diagrams is isomorphic
 to the bialgebra of 3-valent diagrams, is an obvious consequence of the 
 fact that $\calA$ is quasi-isomorphic to $\calD$: the lower line 
 homology of the dual to $\calA$ is the bialgebra of chord diagrams and 
 the lower line
 homology of the dual to $\calD$ is the bialgebra of 3-valent diagrams modulo $STU$, 
 $AS$, and 
 $IHX$ relations.
 
 In Section~\ref{s9} we define a new series of graph-complexes $\{P_n\}_{n\geq 0}$ satisfying
 $
 H^*(P_n)=L_n=Hom(\pi_*(C^n),\Q).
 $
 We show that the totalization complex
 $
 \calP=\Tot P_\bullet
 $
 is quasi-isomorphic to $\Tot L_\bullet$ and therefore
 $
 H^*(\calP)=Hom(\pi_*(\overline{Emb}),\Q).
 $

\section{Cohomology graph-complex for configuration and knot spaces}\label{s8}

Our definition of the space $D_n$ of diagrams is very close to that of~\cite{Kon-OMDQ,LambVol-FLBO}.

A {\it diagram} $\Gamma$ on $n$ external and $q$ internal vertices is any graph with $n$ {\it external} vertices (lying
on the line $\R^1$ and labeled consequently $1,2,\ldots,n$) and $q$ (non-labeled) internal vertices, and some number of oriented segments 
connecting them. Those segments that connect two external vertices are called {\it chords} and all others are {\it edges}.

The {\it orientation set} of a diagram is the union of the set of internal vertices (such elements are considered to be of degree $-d$)
and the set of edges (such elements are of degree $d-1$). An {\it orientation} of a diagram is any ordering of its orientation set.
The {\it degree} of a diagram is the total degree of the elements from the orientation set.

\begin{definition}\label{d81}
A diagram is called {\it admissible} if

(1) it does not contain an internal vertex of valence $\leq 2$;

(2) it contains neither edges nor chords connecting a vertex to itself (no loops);

(3) every internal vertex is connected by a path to an external one.
\end{definition}

\begin{remark}\label{r82}
The distinction between our definition and the one given in~\cite{Kon-OMDQ,LambVol-FLBO} is that we do permit multiple edges and multiple chords. This will be important 
for Theorem~\ref{t93}. This difference is essential only if $d$ is odd, for even $d$ graphs with multiple edges/chords cancel out  by the orientation
relation, see below.
\end{remark}

\begin{definition}\label{d83}
The space $D_n$ is defined as the $\Q$-vector space spanned by the admissible diagrams $\Gamma$ with $n$ external vertices, modulo the relations

(1) if $\Gamma_1$ and $\Gamma_2$ differ only by an orientation of an edge, then
$$
\Gamma_1=(-1)^d\Gamma_2;
$$

(2) if $\Gamma_1$ and $\Gamma_2$ differ only by a permutation of the orientation set, then
$$
\Gamma_1=\pm \Gamma_2
$$
where the sign is the Koszul sign of permutation (taking into account the degrees of elements).
\end{definition}

$D_0$ is defined to be $\Q$ being spanned by the empty diagram.

The differential in $D_n$ is defined as the sum of contractions of edges.

\vspace{0.3cm}

\begin{figure}[h]
\psfrag{+}[0][0][1][0]{$\pm$}
\psfrag{d}[0][0][1][0]{$d_{D_3}$}
\includegraphics[width=14cm]{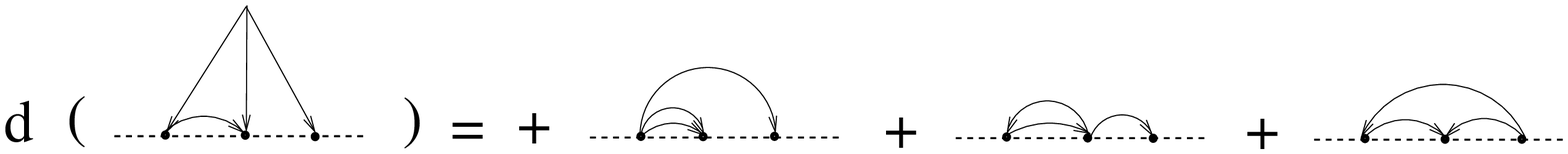}
\caption{The differential in $D_3$.}\label{fig6}
\end{figure}


For the signs convention, see~\cite{Kon-OMDQ,LambVol-FLBO}.

The multiplication in $D_n$ is defined by superimposing:

\vspace{0.3cm}

\begin{figure}[h]
\includegraphics[width=10cm]{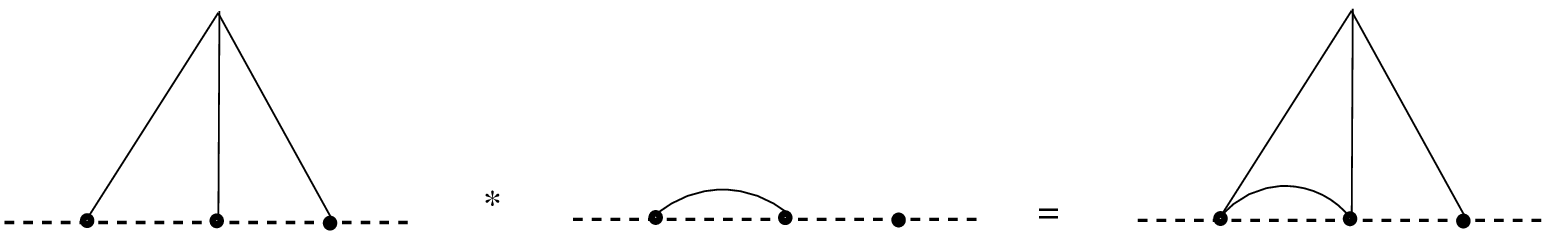}
\caption{The product in $D_3$.}\label{fig7}
\end{figure}


\begin{proposition}\label{p84}
{\rm \cite{Kon-OMDQ,LambVol-FLBO}} Complexes $D_n$, $n\geq 0$,  with multiplication as above are commutative $dg$-algebras.
\end{proposition}

\begin{lemma-def}\label{ld85}
The morphisms
$$
\bar I_n\colon D_n\to A_n
\eqno(\numb)\label{eq81}
$$
which send any diagram having internal vertices to zero, and all others to the corresponding monomials in $A_n$ (see Section~\ref{s3}),
are quasi-isomorphisms of commutative $dg$-algebras.
\end{lemma-def}

\begin{proof}
Our complexes are slightly different from those used in~\cite{Kon-OMDQ,LambVol-FLBO}, but the proof is the same, see~\cite{LambVol-FLBO}.
\end{proof}

The complexes $\{D_n\}_{n\geq 0}$ form a cooperad in $\CDGA$. For the definition of structure maps see~\cite{Kon-OMDQ,LambVol-FLBO}.
This cooperad is endowed with a morphism to the cooperad $\COASS$:
$$
D_\bullet\to\COASS.
$$
Any non-trivial diagram is sent to zero and the trivial diagram with $n$ external vertices --- to  $1\in\Q=\COASS(n)$. This endows $D_\bullet$ with a structure of a simplicial commutative 
$dg$-algebra, see Section~\ref{s2}.
The simplicial structure of $D_\bullet$ is completely analogous to that of $A_\bullet$, see Figures~\ref{fig3}-\ref{fig4}. For $\Gamma\in D_n$, the face map $d_0$ removes the vertex 1 if it was of valence 0,
otherwise $d_0(\Gamma)=0$. Face $d_i(\Gamma)$, $i=1\ldots n-1$, is obtained by contracting the segment $[i,i+1]$ of $\R^1$. And finally $d_n$
removes the last point $n$ (if it was of valence 0, otherwise $d_n(\Gamma)=0$).
The degeneracy $s_i$, $i=1\ldots n+1$, is defined as insertion of a new external point between $i-1$ and $i$.

The normalized part $ND_\bullet$ of $D_\bullet$ is spanned by the diagrams whose all external vertices are of positive valence.
We will define a graph-complex $\calD$ as the totalization of $D_\bullet$.

\begin{theorem}\label{t88}
The complex $\calD=\Tot D_\bullet$ is quasi-isomorphic to $\calA=\Tot A_\bullet$ and therefore
the homology of $\calD$ is  the rational cohomology of $\overline{Emb}$:
$$
H^*(\calD)=H^*(\overline{Emb}).
$$
\end{theorem}

\begin{proof}
The map \eqref{eq81} is a quasi-isomorphism of simplicial commutative $dg$-algebras, which induces a quasi-isomorphism of totalizations:
$$
\calD=\Tot(D_\bullet)\stackrel{\simeq}{\longrightarrow}\Tot(A_\bullet)=\calA.
$$
But $H^*(\calA)=E_2^{*,*}(C^\bullet)=H^*(\overline{Emb})$.
\end{proof}

\begin{figure}[h]
\psfrag{+}[0][0][1][0]{$\pm$}
\psfrag{d}[0][0][1][0]{$d_{\calD}$}
\includegraphics[width=14cm]{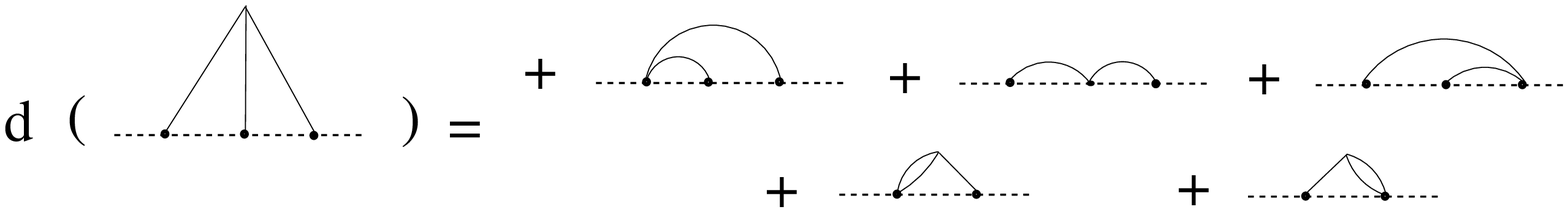}
\caption{The differential in $\calD$.}\label{fig10}
\end{figure}

Since $d_0$ and $d_n$ act always as zero on $ND_n$, the differential in $\calD$ is the sum of contractions of edges and of line segments of $\R^1$, see Figure~\ref{fig10}.

\vspace{0.3cm}

In $\calD=\Tot D_\bullet$ a degree of a graph $\Gamma\in D_n$ is desuspended by $n$. Geometrically we add to the orientation set of $\Gamma$ $n$
elements of degree $-1$ that correspond to the external vertices of $\Gamma$.
The product in $\calD$, which is defined via the Eilenberg-MacLane map, acts as a shuffle of external points. For each summand the ordering of its orientation set is obtained by concatenation.

\vspace{0.3cm}

\begin{figure}[h]
\includegraphics[width=14cm]{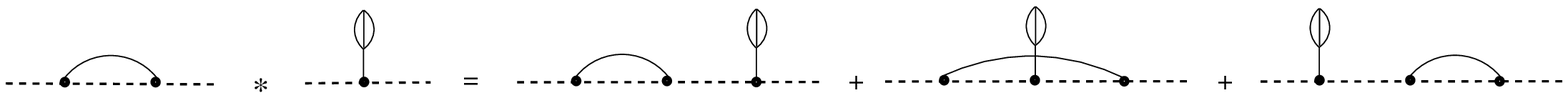}
\caption{The product in $\calD$.}\label{fig11}
\end{figure}

The coproduct in $\calD$ is the coconcatenation.

\vspace{0.3cm}

\begin{figure}[h]
\psfrag{x}[0][0][1][0]{$\otimes$}
\psfrag{D}[0][0][1][0]{$\Delta$}
\includegraphics[width=14cm]{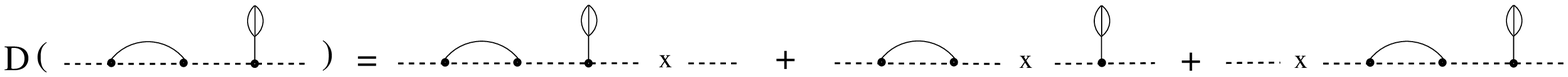}
\caption{The coproduct in $\calD$.}\label{fig12}
\end{figure}

\section{Cohomotopy graph-complex for configuration and knot spaces}\label{s9}
A non-trivial graph $\Gamma\in D_n$ is called {\it non-decomposable} if it can not be represented as a product $\Gamma=\Gamma_1\cdot\Gamma_2$
of two non-trivial graphs $\Gamma_1,\Gamma_2\in D_n$.\footnote{We consider the inner product of $D_n$, see Figure~\ref{fig7}.} The space spanned 
by non-decomposable graphs will be denoted by $P_n$.

\begin{remark}\label{r91}
In other words a graph is non-decomposable if it is connected (and non-empty) when we remove from it all the
external vertices with their little neighborhoods.
\end{remark}

Notice however that a non-decomposable graph might be disconnected: together with its main connected part it can have a number of singletons --- external vertices of valence~0.

The complex $P_n$ is a quotient-complex of $D_n$. It is easy to see that $P_n=D_n/D_n^{\,2}$.

\begin{proposition}\label{p92}
$D_n$ is a polynomial algebra whose space of generators is $P_n$.
\end{proposition}

\begin{proof}
Obvious. It is here where it is important that we permit multiple edges/chords.
\end{proof}

\begin{theorem}\label{t93}
The homology of $P_n$ is the rational cohomotopy of the configuration space:
$$
H^*(P_n)=L_n=Mor(\pi_*(C^n),\Q).
$$
\end{theorem}

\begin{proof}
We have the quasi-isomorphisms
$$
\xymatrix{
P_n=D_n/D_n^{\,2}&\calL(D_n)\ar@{->>}[l]^-\alpha_-\simeq\ar@{->>}[r]^-\simeq&\calL(A_n)&L_n
\ar@{_{(}->}[l]_-\simeq.
}
\eqno(\numb)\label{eq91}
$$
The first arrow $\alpha$ is a quasi-isomorphism by Proposition~\ref{p52}, the second one is induced by the quasi-isomorphism~\ref{eq81}, the last one 
is due to the Koszul property.
\end{proof}

It follows from Remark~\ref{r91} that $P_\bullet$ is a simplicial subspace of $D_\bullet$ ---
the simplicial structure maps preserve $P_\bullet$. Its noramalized part $NP_\bullet$ is spanned
by the connected non-decomposable diagrams, \ie by the diagrams without singletons. 
Now define complex $\calP$ as the totalization of $P_\bullet$. Obviously, $\calP$ is a quotient-complex of $\calD$. 

\vspace{0.3cm}

\begin{figure}[h]
\psfrag{+}[0][0][1][0]{$\pm$}
\psfrag{d}[0][0][1][0]{$d_{\calP}$}
\includegraphics[width=11cm]{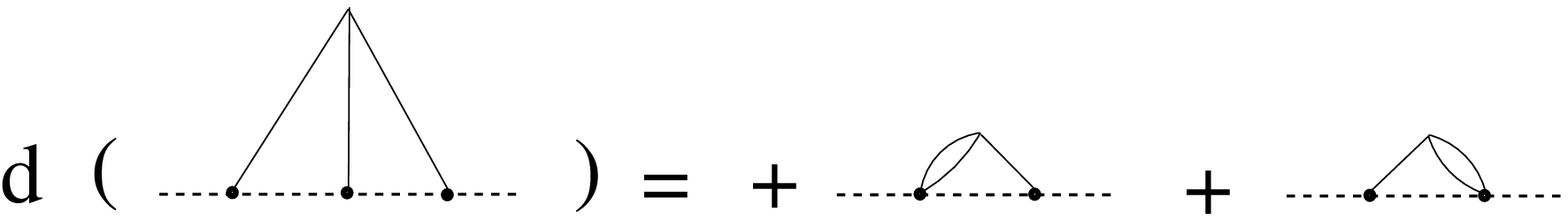}
\caption{The differential in $\calP$.}\label{fig13}
\end{figure}

\begin{theorem}\label{t94}
The complex $\calP=\Tot P_\bullet$ is quasi-isomorphic to $\Tot L_\bullet$ and therefore the homology of $\calP$ is the rational cohomotopy of $\overline{Emb}$:
$$
H^*(\calP)=Mor(\pi_*(\overline{Emb}),\Q).
$$
\end{theorem}

\begin{proof}
Diagram \eqref{eq91} is a sequence of quasi-isomorphisms of simplicial $dg$-spaces. Passing to totalization one gets the result.
\end{proof}

\begin{remark}
$P_\bullet$ is a simplicial $L_\infty$-coalgebra. Indeed, given $B$ is a polynomial $dg$-algebra,
any section $B/B^2\hookrightarrow B$ of the projection $B\twoheadrightarrow B/B^2$ defines an $L_\infty$-coalgebra structure on $B/B^2$. We have natural inclusions 
$$
D_n/D_n^{\,2}=P_n\hookrightarrow D_n.
\eqno(\numb)\label{eq92}
$$
Since~\eqref{eq92} is a map of simplicial vector spaces $P_\bullet\to D_\bullet$, the $L_\infty$-coalgebra operations on $P_n$, $n\geq 0$, commute with the simplicial structure
maps.
\end{remark}

We finish by giving some examples of cycles in $\calP$. Obviously the diagram 
$$
\includegraphics[width=3cm]{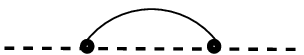}
$$
is the first non-trivial cycle (for both even and odd $d$). Its degree is $d-3$. It can be easily seen 
that the sum of diagrams (taken with appropiate signs)
$$
\includegraphics[width=7cm]{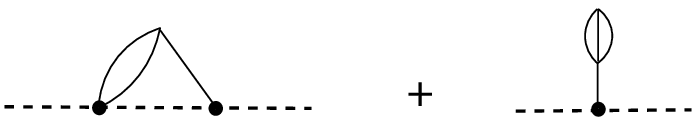}
$$
is a non-trivial cycle in case when $d$ is odd (and therefore multiple edges are possible). The degree of this cycle is $2d-5$.
Recall that $\overline{Emb}$ is homotopy equivalent to $Emb\times \Omega^2S^{d-1}$. 
The above cycles descibe the rational cohomotopy coming
from the second factor $\Omega^2S^{d-1}$.

The first non-trivial cohomotopy coming from the first factor $Emb$ is of degree $2d-6$. In case of even $d$ it is given by the diagram:
$$
\includegraphics[width=2.3cm]{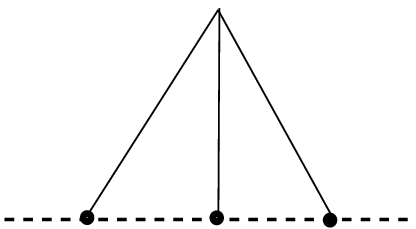}
$$
For odd $d$ it is given by the sum:
$$
\includegraphics[width=5cm]{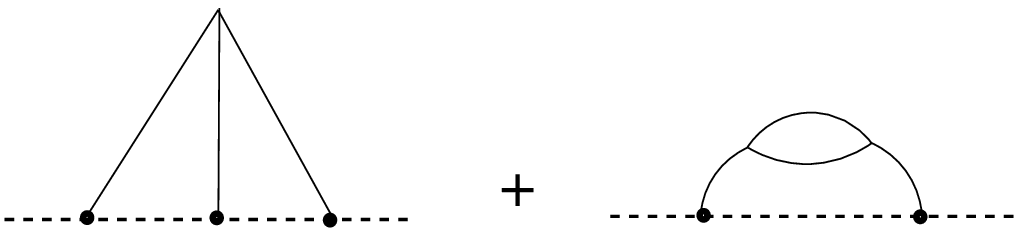}
$$

\subsection*{Acknowledgements}
The second author thanks D.~Sinha, J.~Conant, and L.~Ionescu  for fruitful conversations. He is also
grateful to the University of Oregon and IHES where this work was completed for hospitality. We thank R.~Budney for an attentive reading and various corrections.

\end{document}